\documentclass[11pt,abstracton]{scrartcl}
\usepackage{booktabs}
\usepackage{epsfig}
\usepackage{graphicx}
\usepackage{amssymb}
\usepackage{amsfonts}
\usepackage{multirow}
\usepackage{algorithmic}
\usepackage{amsmath}
\usepackage{amsthm}
\usepackage{mathrsfs}
\usepackage[linesnumbered,ruled,vlined]{algorithm2e}
\usepackage{mathtools}
\usepackage[dvipsnames,table]{xcolor}
\usepackage{subfigure}
\usepackage{tikz}
\usepackage{paralist}
\usepackage{color}
\usepackage{tabu}
\usepackage{rotating}
\usepackage{pifont}
\usepackage{cite}
\usepackage{threeparttable}
\usepackage{moresize}
\usepackage{cancel}
\usepackage{pgf-pie}
\usepackage{adjustbox}
\usepackage[margin=2cm]{geometry}
\usepackage[noblocks,affil-sl]{authblk}
\usepackage{url}
\usepackage{blkarray}

\usepackage[colorlinks=true,linkcolor=black,anchorcolor=black,citecolor=black,filecolor=black,menucolor=black,runcolor=black,urlcolor=black]{hyperref}

\renewcommand{\vec}{\mathbf}

\title{Generating Large-scale Dynamic Optimization Problem Instances Using the Generalized Moving Peaks Benchmark}
\author[1]{Mohammad~Nabi~Omidvar\thanks{m.n.omidvar@leeds.ac.uk}}
\author[2]{Danial~Yazdani\thanks{danial.yazdani@gmail.com,yazdani@sustech.edu.cn}}
\author[3]{J\"{u}rgen~Branke\thanks{Juergen.Branke@wbs.ac.uk}}
\author[4]{ Xiaodong~Li\thanks{xiaodong.li@rmit.edu.au}}
\author[5]{Shengxiang~Yang\thanks{syang@dmu.ac.uk}}
\author[2,7]{Xin~Yao\thanks{xiny@sustc.edu.cn}}
\affil[1]{\normalfont School of Computing, University of Leeds, and Leeds University Business School, Leeds, United Kingdom.}
\affil[2]{\normalfont Guangdong Provincial Key Laboratory of Brain-inspired Intelligent Computation, Department of Computer Science and Engineering, Southern University of Science and Technology, Shenzhen 518055, China.}
\affil[3]{\normalfont Operational Research and Management Sciences Group, Warwick Business School, University of Warwick, Coventry CV4 7AL, United Kingdom.}
\affil[4]{\normalfont School of Science (Computer Science), RMIT University, GPO Box 2476, Melbourne, 3001, Australia.}
\affil[5]{\normalfont School of Computer Science and Informatics, De Montfort University, Leicester, United Kingdom.}
\affil[6]{\normalfont Guangdong Provincial Key Laboratory of Brain-inspired Intelligent Computation, Department of Computer Science and Engineering, Southern University of Science and Technology, Shenzhen 518055, China.}
\affil[7]{Center of Excellence for Research in Computational Intelligence and Applications (CERCIA), School of Computer Science, University of Birmingham, Birmingham B15 2TT, United Kingdom.}
\date{\large July 2021}

\providecommand{\keywords}[1]
{
  \small	
  \textbf{\textit{Keywords---}} #1
}

\begin{document}

\maketitle

\begin{abstract}
This document describes the generalized moving peaks benchmark (GMPB)~\cite{yazdani2020benchmarking} and how it can be used to generate problem instances for continuous large-scale dynamic optimization problems. It presents a set 15 benchmark problems, the relevant source code, and a performance indicator, designed for comparative studies and competitions in large-scale dynamic optimization. Although its primary purpose is to provide a coherent basis for running competitions, its generality allows the interested reader to use this document as a guide to design customized problem instances to investigate issues beyond the scope of the presented benchmark suite. To this end, we explain the modular structure of the GMPB and how its constituents can be assembled to form problem instances with a variety of controllable characteristics ranging from unimodal to highly multimodal, symmetric to highly asymmetric, smooth to highly irregular, and various degrees of variable interaction and ill-conditioning.
\end{abstract}


\section{Introduction}

Change is an inescapable aspect of natural and artificial systems, and adaptation is central to their resilience.
Optimization problems are no exception to this maxim.
Indeed, viability of businesses and their operational success depend heavily on their effectiveness in responding to a change in the myriad of optimization problems they entail. 
For an optimization problem, this boils down to the efficiency of an algorithm to find and maintain a sequence of quality solutions to an ever changing problem.

Ubiquity of dynamic optimization problems (DOPs)~\cite{yazdani2021DOPsurveyPartA} demands extensive research into design and development of algorithms capable of dealing with various types of change~\cite{nguyen2012evolutionary}.
These are often attributed to a change in the objective function, its number of decision variables, or constraints.
Despite the large body of literature on dynamic optimization problems and algorithms, little attention has been given to their scalability~\cite{yazdani2019scaling,yazdani2021DOPsurveyPartB}. 
Indeed, the number of dimensions of a typical DOP studied in the literature rarely exceeds twenty. 

Motivated by rapid technological advancements, large-scale optimization has gained popularity in recent years~\cite{mahdavi2015metaheuristics}.
However, the exponential growth in the size of the search space, with respect to an increase in the number of decision variables, has made large-scale optimization a challenging task.
For DOPs, however, the challenge is twofold.
For such problems, not only should an algorithm be capable of finding the global optimum in the vastness of the search space but it should also be able to track it over time.
For multi-modal DOPs, where several optima have the potential to turn into the global optimum after environmental changes, the cost of tracking multiple moving optima also adds to the complexity.

Real-world large-scale optimization problems often exhibit a modular structure with nonuniform imbalance among the contribution of its constituent parts to the objective value~\cite{omidvar2015designing,omidvar2017dg2}. 
The modularity is caused by the interaction structure of the decision variables resulting in a wide range of structures from fully separable functions to fully nonseparable ones. 
Most problems exhibit some degree of sparsity in their interaction structure, which can be exploited by optimization algorithms. 
The imbalance property can be caused as a by-product of modularity or due to the heterogeneous nature of the input variables and their domains. 

Generalized moving peaks benchmark (GMPB)~\cite{yazdani2020benchmarking} is capable of generating problem instances with a variety of characteristics that can range from fully non-separable to fully separable structure, from homogeneous to highly heterogeneous sub-functions, and from balanced to highly imbalanced sub-functions.
Each sub-function generated by GMPB is constructed by assembling several components.
In most benchmark generators in the filed of DOPs, these components are unimodal, smooth, symmetric, fully separable, and easy-to-optimize \emph{peaks}.
However, GMPB is capable of generating components with a variety of properties that can range from unimodal to highly multimodal, from symmetric to highly asymmetric, from smooth/regular to highly irregular, with different variable interaction degrees, and from low condition number to highly ill-conditioned.

\section{Generalized Moving Peaks Benchmark~\cite{yazdani2020benchmarking,yazdani2021generalized}}
\label{sec:GMPB}

GMPB's main function is constructed by assembling several sub-functions as follows:
\begin{align}
\label{eq:wCMPB}
F^{(t)}(\vec{x}) = d^{-1}\sum_{i=1}^n \omega_i d_i f_i^{(t)}(\vec{x}),
\end{align}
where $t$ shows the current environment number, $f_i^{(t)}$ is the $i$th sub-function in the $t$th environment, $n$ is the number of sub-functions, $d$ is the dimension of the main function, $d_i$ is the dimension of the sub-function $f_i$, and $\omega_i$ controls the contribution of sub-function $f_i$ for generating imbalance property.
The baseline function that generates each sub-function $f_i$ in GMPB is:
\begin{align}
\label{eq:irGMPB}
 f_i^{(t)}(\vec{x}_i)= \max_{k\in\{1,\dots,m_i\}}\left\{ h_k^{(t)} - \sqrt{\mathbb{T}\left(\left(\vec{x}_i-\vec{c}_{k,i}^{(t)}\right)^\top{\mathbf{R}_{k,i}^{(t)}}^\top,k,i\right)  \mathbf{W}_{k,i}^{(t)}  \mathbb{T}\left(\mathbf{R}_{k,i}^{(t)}\left(\vec{x}_i-\vec{c}_{k,i}^{(t)}\right),k,i\right)} \right\},
\end{align}
where $\mathbb{T}(\vec{y},k,i):\mathbb{R}^{d_i}\mapsto\mathbb{R}^{d_i}$ is calculated as~\cite{Hansen2010real}: 
\begin{align}
\label{eq:ir}
    \mathbb{T}\left(y_{j},k,i\right)=
    \begin{dcases}
    \exp{\left(\log(y_{j})+\tau^{(t)}_{k,i}\left(\sin{(\eta_{k,i,1}^{(t)}\log(y_{j}))}+\sin{(\eta_{k,i,2}^{(t)}\log(y_{j}))} \right)\right)} & \text{if   } y_{j}>0 \\
    0 & \text{if   }y_{j}=0\\
    -\exp{\left(\log(|y_{j}|)+\tau^{(t)}_{k,i}\left(\sin{(\eta_{k,i,3}^{(t)}\log(|y_{j}|))}+\sin{(\eta_{k,i,4}^{(t)}\log(|y_{j}|))} \right)\right)} & \text{if   } y_{j}<0  
\end{dcases}
\end{align}
where $\vec{x}_i$ is a subset of decision variables ($d_i$-dimensional) of the solution $\vec{x}$, which belongs to the $i$th sub-function,  $m_i$ is the number of components in $f_i$, $\mathbf{R}_{k,i}^{(t)}$ is the rotation matrix of the $k$th component  of the $i$th sub-function in the $t$th environment, $\mathbf{W}_{k,i}^{(t)}$ is a $d_i \times d_i$ diagonal matrix whose elements determine the width of the $k$th component in different dimensions,  
$y_{j}$ is $j$th element of $\mathbf{y}$, and $\eta_{k,i,l\in\{1,2,3,4\}}^{(t)}$ and $\tau^{(t)}_{k,i}$  are irregularity parameters of the $k$th component of the $i$th sub-function.

For each component $k$ of the $i$th sub-function, the rotation matrix $\mathbf{R}_{k,i}$ is obtained by rotating the projection of $\vec{x}_i$ onto all $x_p$-$x_q$ planes by a given angle $\theta_{k,i}$. 
The total number of unique planes which will be rotated is ${d_i \choose 2} = \frac{d_i(d_i-1)}{2}$. 
For rotating each  $x_p$-$x_q$ plane by a certain angle ($\theta_{k,i}$), a \emph{Givens} rotation matrix, $\mathbf{G}_{(p,q)}$, must be constructed.
To do this, $\mathbf{G}_{(p,q)}$ is first initialized to an identity matrix $\mathbf{I}_{d_i \times d_i}$; then, four elements of $\mathbf{G}_{(p,q)}$ are altered as follows: 
\newcommand{\pp}{\cos\left(\theta_{k,i}^{(t)}\right)}
\newcommand{\qq}{\cos\left(\theta_{k,i}^{(t)}\right)}
\newcommand{\pq}{ -\sin\left(\theta_{k,i}^{(t)}\right)}
\newcommand{\qp}{ \sin\left(\theta_{k,i}^{(t)}\right)}
\begin{equation}
\mathbf{G}_{(p,q)} \; = \;
\begin{blockarray}{cccccccc}
       &        &     p     &        &      q        &        &        &   \\
    \begin{block}{(ccccccc)c}
1      & \cdots & 0           & \cdots & 0              & \cdots & 0      &   \\
\vdots & \ddots & \vdots      &        & \vdots         &        & \vdots &   \\
0      & \cdots & \pp         & \cdots & \pq            & \cdots & 0      & p \\
\vdots &        & \vdots      & \ddots & \vdots         &        & \vdots &   \\
0      & \cdots & \qp         & \cdots & \qq            & \cdots & 0      & q \\
\vdots &        & \vdots      &        & \vdots         & \ddots & \vdots &   \\
0      & \cdots & 0           & \cdots & 0              & \cdots & 1      &   \\
\end{block}
\end{blockarray}
     \end{equation}
Thus, the rotation matrix, $\mathbf{R}_{k,i}$, in the $t$th environment is calculated as follows:
\begin{align}
\label{eq:allrotation}
 \mathbf{R}_{k,i}^{(t)}= \prod_{(p,q)\in \mathcal{P}} \mathbf{G}_{(p,q)}   \mathbf{R}_{k,i}^{(t-1)},
\end{align}
where $\mathcal{P}$ contains all unique pairs of dimensions defining all possible planes in a $d_i$-dimensional space. 
The order of the multiplications of the Givens rotation matrices is random.
The reason behind using~\eqref{eq:allrotation} for calculating $\mathbf{R}$ is that we aim to have control on the rotation matrix based on an angle severity $\tilde{\theta}_i$. 
Note that the initial $\mathbf{R}^{(0)}_{k,i}$ for problem instances with rotation property is obtained by using  the Gram-Schmidt orthogonalization method on a matrix with normally distributed entries. 

For each component $k$ in the $i$th sub-function, the height, width vector, center, angle, and irregularity parameters change from one environment to the next according to the following update rules:
\begin{align}
\vec{c}_{i,k}^{(t+1)}&=\vec{c}_{i,k}^{(t)} + \tilde{s}_i \frac{\vec{r}}{\|\vec{r}\|},  \label{eq:center} \\
h_{i,k}^{(t +1)}&=h_{i,k}^{(t)}  + \tilde{h}_i \, \mathcal{N}(0,1), \label{eq:height} \\
w_{i,k,j}^{(t +1)}&=w_{i,k,j}^{(t)}  + \tilde{w}_i \, \mathcal{N}(0,1), j\in\{1,2,\cdots,d_i\},\label{eq:width} \\
\theta_{i,k}^{(t +1)}&= \theta_{i,k}^{(t)}+\tilde{\theta}_i \, \mathcal{N}(0,1),\label{eq:angle}\\
\eta_{i,k,l}^{(t +1)}&= \eta_{i,k,l}^{(t)}+\tilde{\eta}_i \, \mathcal{N}(0,1), l\in\{1,2,3,4\},\label{eq:tau}\\
\tau_{i,k}^{(t +1)}&= \tau_{i,k}^{(t)}+\tilde{\tau}_i \, \mathcal{N}(0,1),\label{eq:eta}
\end{align}
where $\mathcal{N}(0,1)$ is a random number drawn from a Gaussian distribution with mean 0 and variance 1, $\mathbf{c}_{i,k}$ shows the vector of center position of the $k$th components of the $i$th sub-function, $\vec{r}$ is a $d_i$-dimensional vector of random numbers generated by $\mathcal{N}(0,1)$,  $\|\vec{r}\|$ is the Euclidean length (i.e., $l_2$-norm) of $\vec{r}$,  $\frac{\vec{r}}{\|\vec{r}\|}$ generates a unit vector with a random direction,  $\tilde{h}_i$, $\tilde{w}_i$, $\tilde{s}_i$, $\tilde{\theta}_i$, $\tilde{\eta}_i$ and $\tilde{\tau}_i$ are height, width, shift, angle, and irregularity parameters' change severity values of components in the $i$h sub-function, respectively, $w_{i,k,j}$\footnote{Not to be confused with $\omega$ that controls the components' imbalance.} shows the width of the $k$th component in the $j$th dimension of the $i$th sub-function, and $h_{i,k}$ and $\theta_{i,k}$ show the height and angle of the $k$th component in the $i$th sub-function, respectively.

Outputs of equations~\eqref{eq:center}~to~\eqref{eq:eta} are bounded as follows: $h_{i,k}\in[h_{\mathrm{min}},h_\mathrm{max}]$, $w_{i,k,j}\in[w_\mathrm{min},w_\mathrm{max}]$, $\vec{c}_{i,k}\in[Lb_i,Ub_i]^{d_i}$, $\tau\in[\tau_\mathrm{min},\tau_\mathrm{max}]$, $\eta_{1,2,3,4}\in[\eta_\mathrm{min},\eta_\mathrm{max}]$, and $\theta_{i,k}\in[\theta_\mathrm{min},\theta_\mathrm{max}]$, where $Lb_i$ and $Ub_i$ are maximum and minimum problem space bounds in the $i$th sub-function. 
For keeping the above mentioned values in their bounds, a \emph{Reflect} method is utilized. Assume $a^{(t+1)}=a^{(t)}+b$ represents one of the equations~\eqref{eq:center}~to~\eqref{eq:eta}. 
The output based on the reflect method is:
\begin{gather}
\label{chap5:eq:survivalTime}
a^{(t+1)}=
\begin{dcases}
a^{(t)}+b & \text{if   }a^{(t)}+b\in[a_\mathrm{min},a_\mathrm{max}] \\
2\times a_\mathrm{min}-a^{(t)}-b & \text{if   }a^{(t)}+b<a_\mathrm{min} \\
2\times a_\mathrm{max}-a^{(t)}-b & \text{if   }a^{(t)}+b>a_\mathrm{max} 
\end{dcases}
\end{gather}

\subsection{Characteristics of the components, sub-functions, and problem instances}

In this section, we describe the main characteristics of the components, sub-functions, and problem instances generated by GMPB.

\subsubsection{Component characteristics}

\begin{figure}[tp!]
\centering
\begin{tabular}{cc}
     \subfigure[{\scriptsize }]{\includegraphics[width=0.45\linewidth]{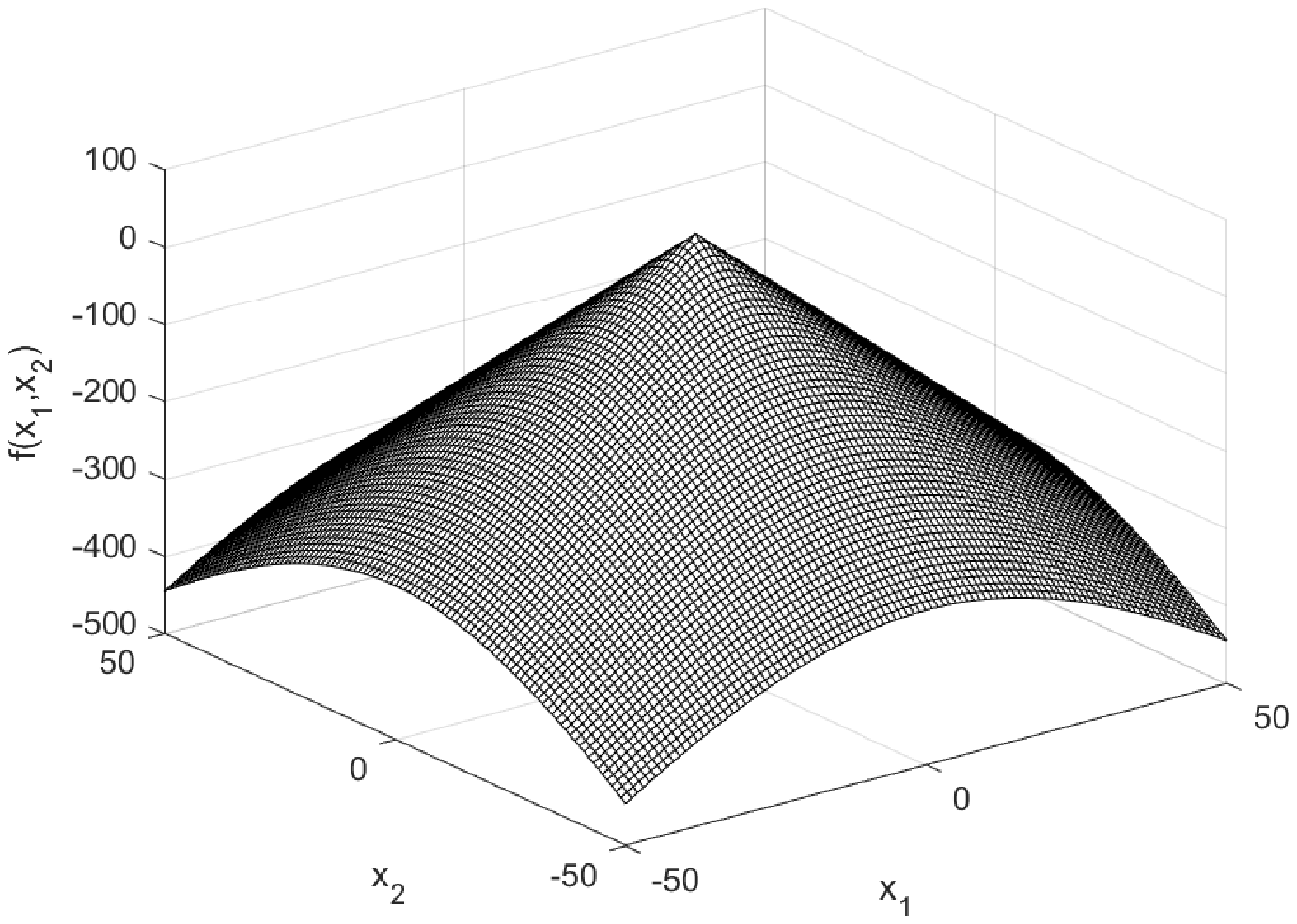}\label{fig:Cmponent:cone:surf}}
&
    \subfigure[{\scriptsize }]{\includegraphics[width=0.45\linewidth]{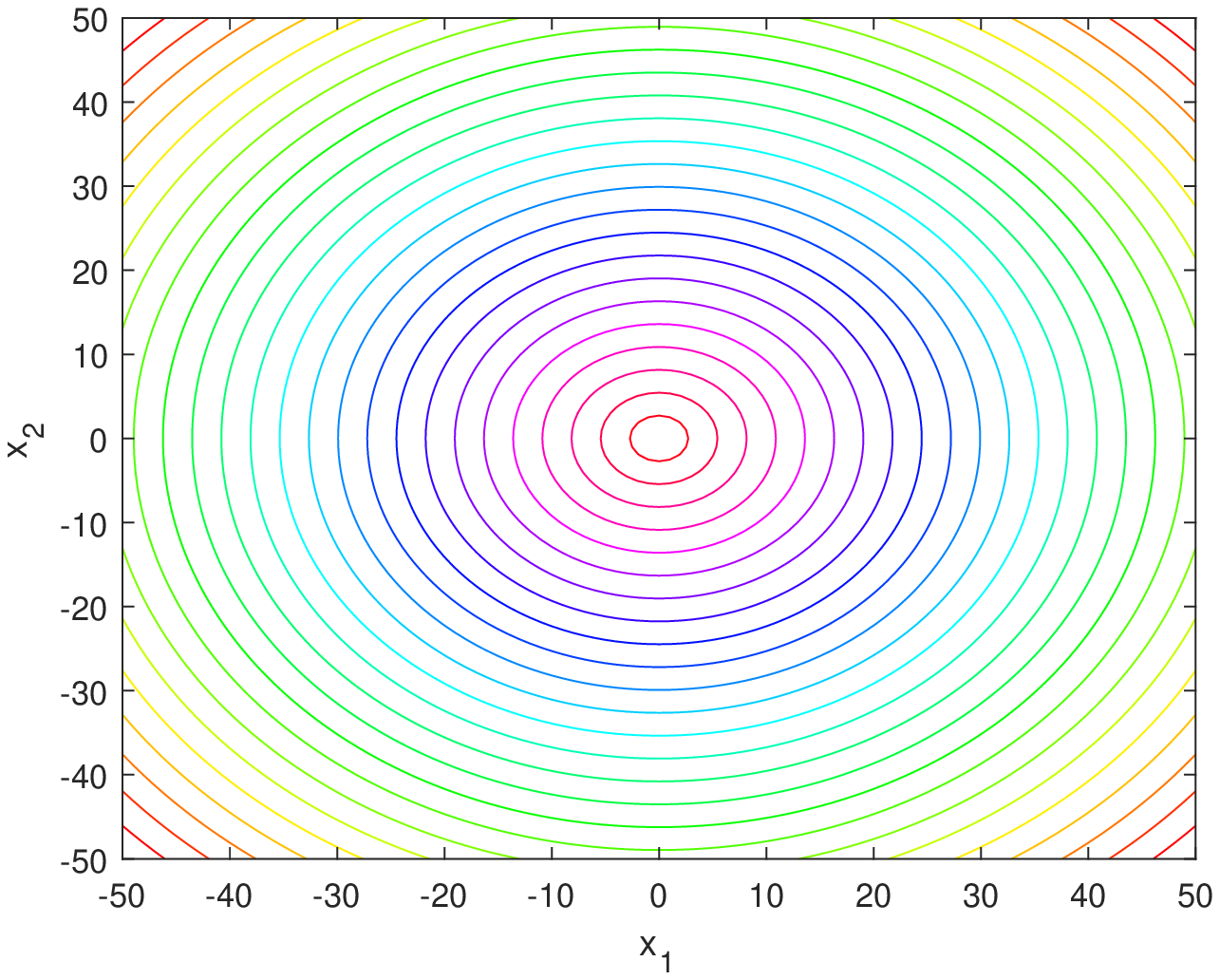}\label{fig:Cmponent:cone:contour}}
\end{tabular}
\caption{A component generated by \eqref{eq:irGMPB} whose width values are identical, $\mathbf{R}=\mathbf{I}$, and $\tau$ and $\eta$ are set to zero.}
\label{fig:Cmponent:cone}
\end{figure}

\begin{figure}[tp!]
\centering
\begin{tabular}{cc}
     \subfigure[{\scriptsize $\tau=0.2$ and $\eta_{1,2,3,4}=[15,15,15,15]$}]{\includegraphics[width=0.45\linewidth]{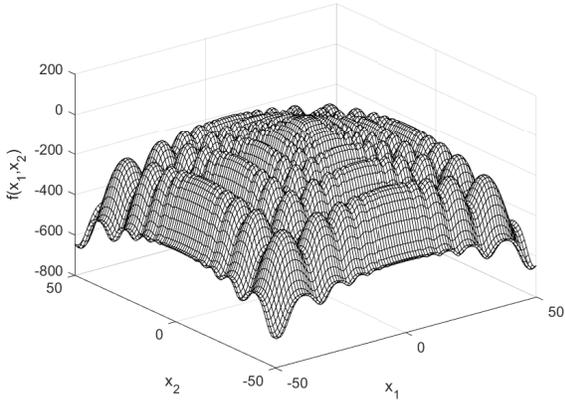}\label{fig:Cmponent:irregular:surf}}
&
    \subfigure[{\scriptsize }]{\includegraphics[width=0.45\linewidth]{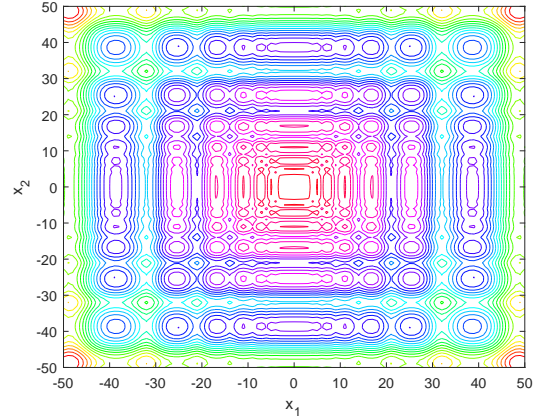}\label{fig:Cmponent:irregular:contour}}
    \\
         \subfigure[{\scriptsize $\tau=0.5$ and $\eta_{1,2,3,4}=[5,5,5,5]$}]{\includegraphics[width=0.45\linewidth]{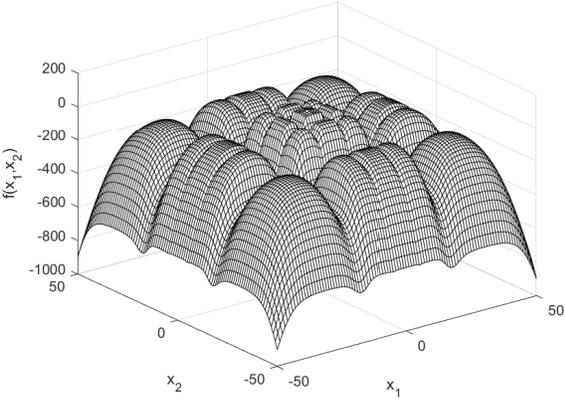}\label{fig:Cmponent:irregular(Challenging):surf}}
&
    \subfigure[{\scriptsize }]{\includegraphics[width=0.45\linewidth]{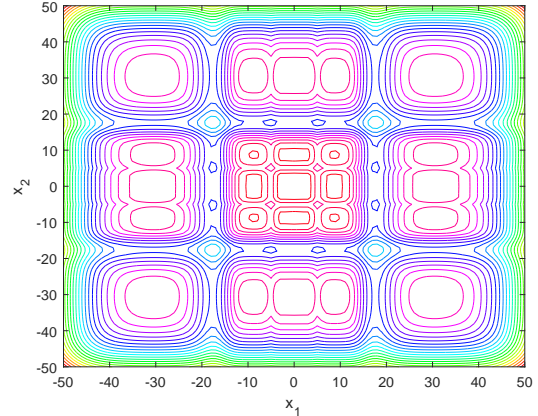}\label{fig:Cmponent:irregular(Challenging):contour}}
\end{tabular}
\caption{Two components generated by \eqref{eq:irGMPB} whose width values are identical and $\mathbf{R}=\mathbf{I}$.}
\label{fig:Cmponent:irregular}
\end{figure}

\begin{figure}[tp!]
\centering
\begin{tabular}{cc}
     \subfigure[{\scriptsize }]{\includegraphics[width=0.45\linewidth]{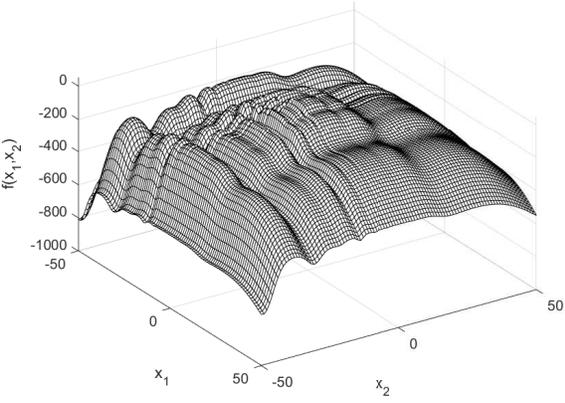}\label{fig:Component_irregular-assymetric:surf}}
&
    \subfigure[{\scriptsize }]{\includegraphics[width=0.45\linewidth]{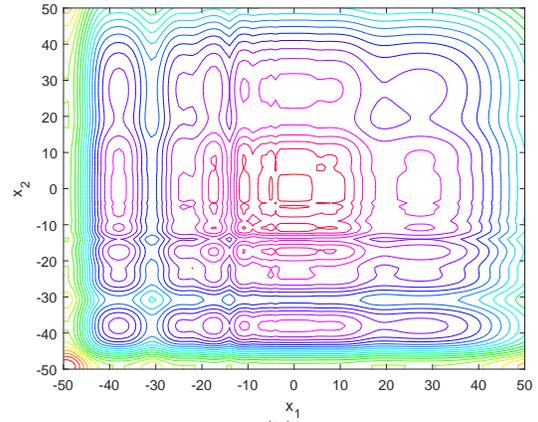}\label{fig:Component_irregular-assymetric:contour}}
\end{tabular}
\caption{An asymmetric component generated by \eqref{eq:irGMPB} whose width values are identical, $\mathbf{R}=\mathbf{I}$, $\tau=0.3$, and $\eta_{1,2,3,4}=[0,5,10,15]$.}
\label{fig:Cmponent:assymetric}
\end{figure}

In the simplest form, \eqref{eq:irGMPB} generates a symmetric, unimodal, smooth/regular, and easy-to-optimize conical peak (see Figure~\ref{fig:Cmponent:cone}).
By setting different values to  $\tau$ and $\eta$, a component generated by GMPB becomes irregular and multimodal.
Figure~\ref{fig:Cmponent:irregular} shows two irregular and multimodal components generated by GMPB with different parameter settings for $\tau$ and $\eta$.
The illustrated component in Figure~\ref{fig:Cmponent:irregular:surf} is easier to optimize in comparison to the one shown in Figure~\ref{fig:Cmponent:irregular(Challenging):surf}.
As can be seen in Figure~\ref{fig:Cmponent:irregular(Challenging):contour}, the local optima on the basin of attraction of the component are very vast which results in premature convergences.
The results shown in the supplementary document of~\cite{yazdani2020benchmarking} indicate that the performance of the optimization algorithms such as the particle swarm optimization (PSO)~\cite{kennedy1995particle} and differential evolution (DE)~\cite{brest2006selfadaptive} deteriorates significantly in such components.
A component is symmetric when all $\eta$ values are identical, and components whose $\eta$ values are set to different values are asymmetric (see Figure~\ref{fig:Cmponent:assymetric}).

\begin{figure}[tp!]
\centering
\begin{tabular}{cc}
     \subfigure[{\scriptsize $\mathbf{w}=[7,7]$ (without ill-conditioning).}]{\includegraphics[width=0.45\linewidth]{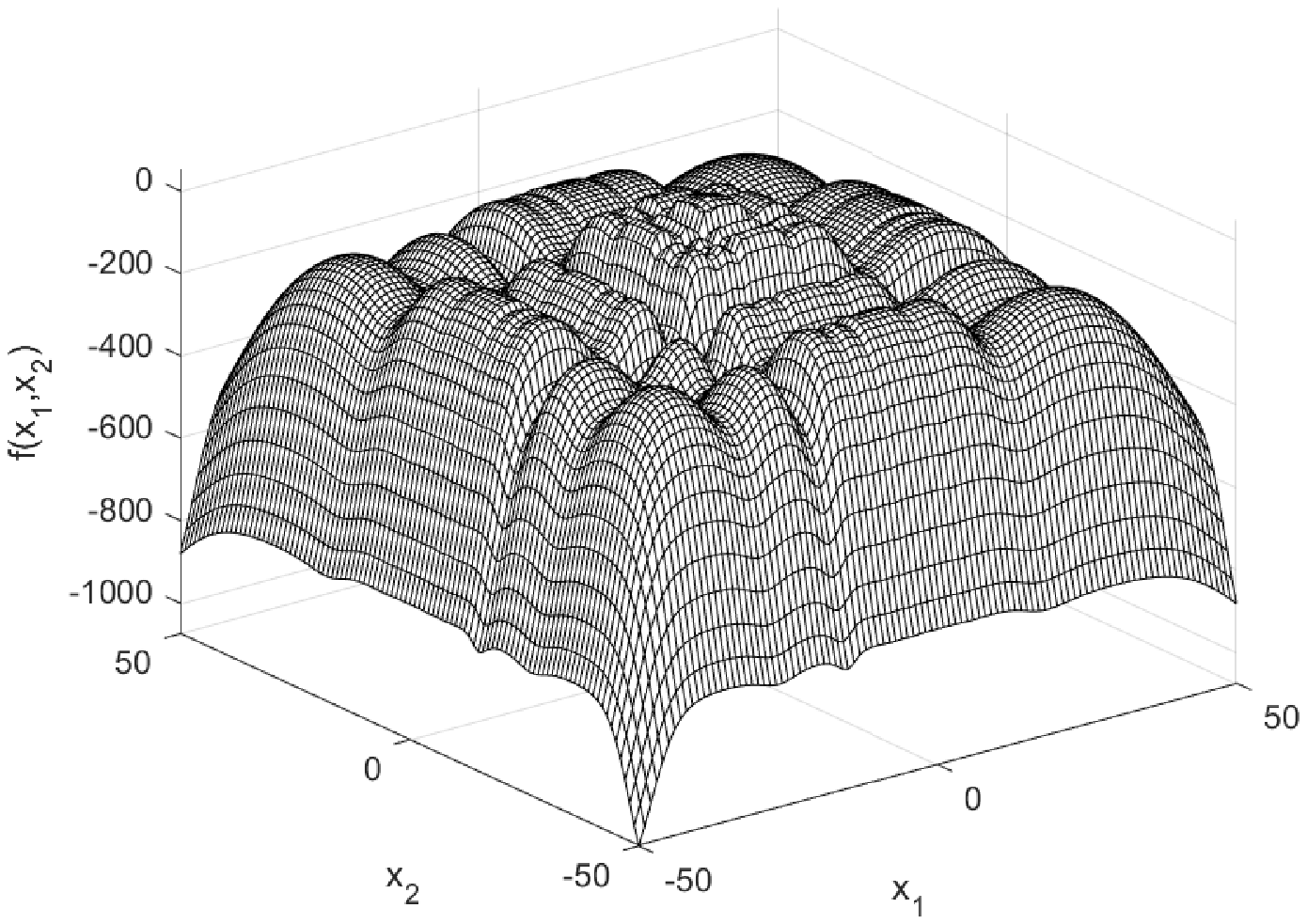}\label{fig:Component_irregular-assymetric-w77-surf}}
&
    \subfigure[{\scriptsize }]{\includegraphics[width=0.45\linewidth]{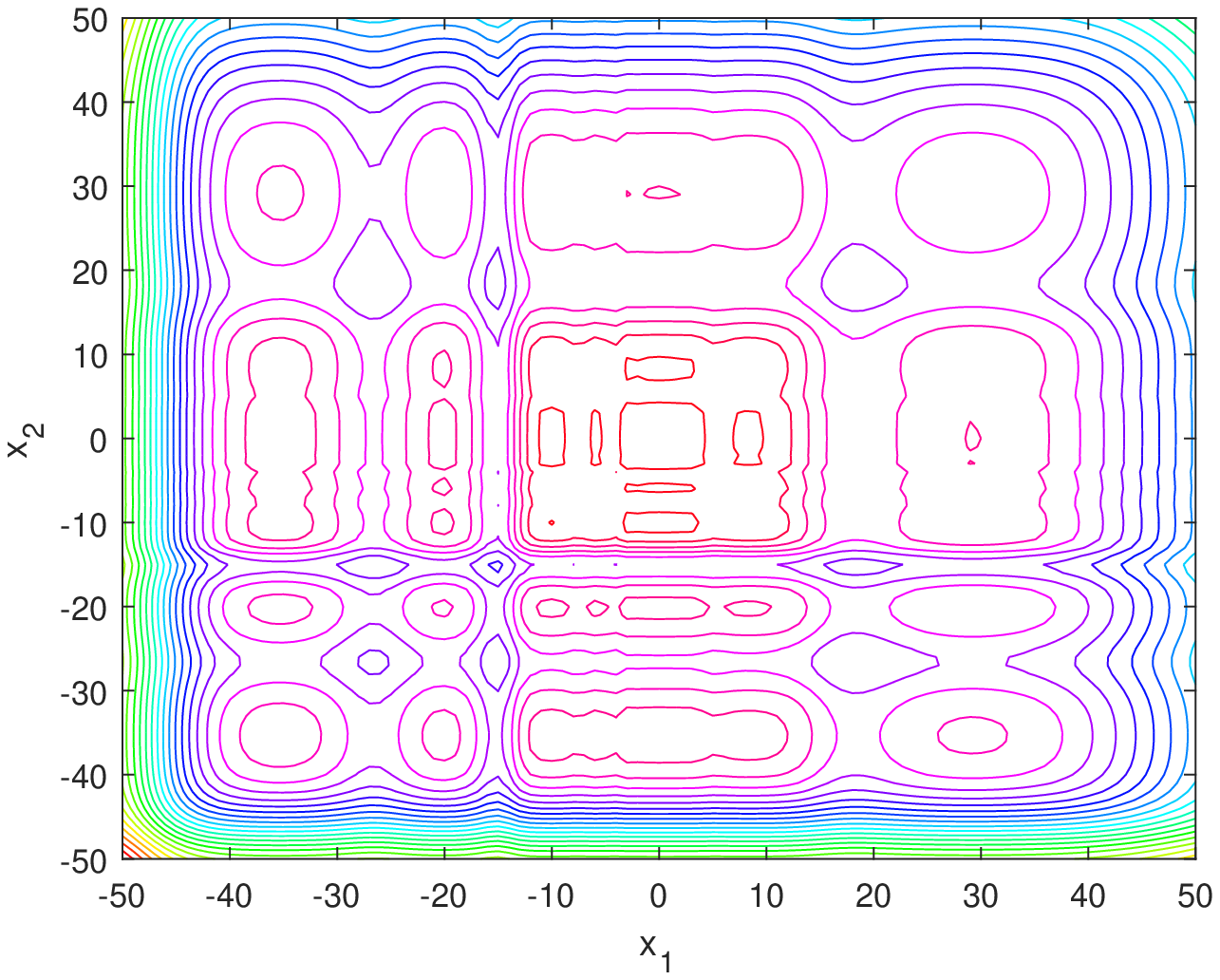}\label{fig:Component_irregular-assymetric-w77-contour}}
    \\
         \subfigure[{\scriptsize $\mathbf{w}=[7,3]$ (ill-conditioned).}]{\includegraphics[width=0.45\linewidth]{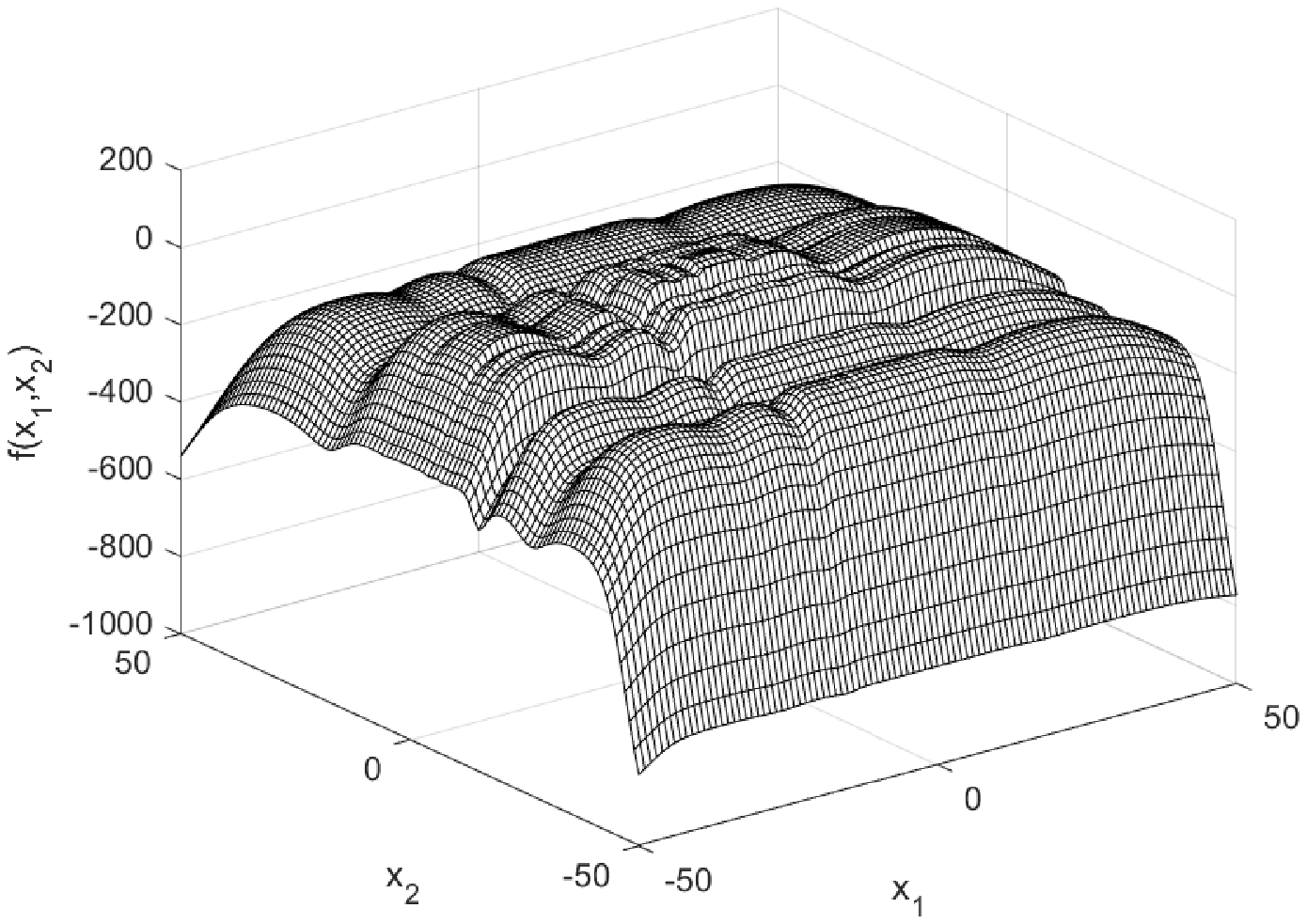}\label{fig:Component_irregular-assymetric-w73-surf}}
&
    \subfigure[{\scriptsize }]{\includegraphics[width=0.45\linewidth]{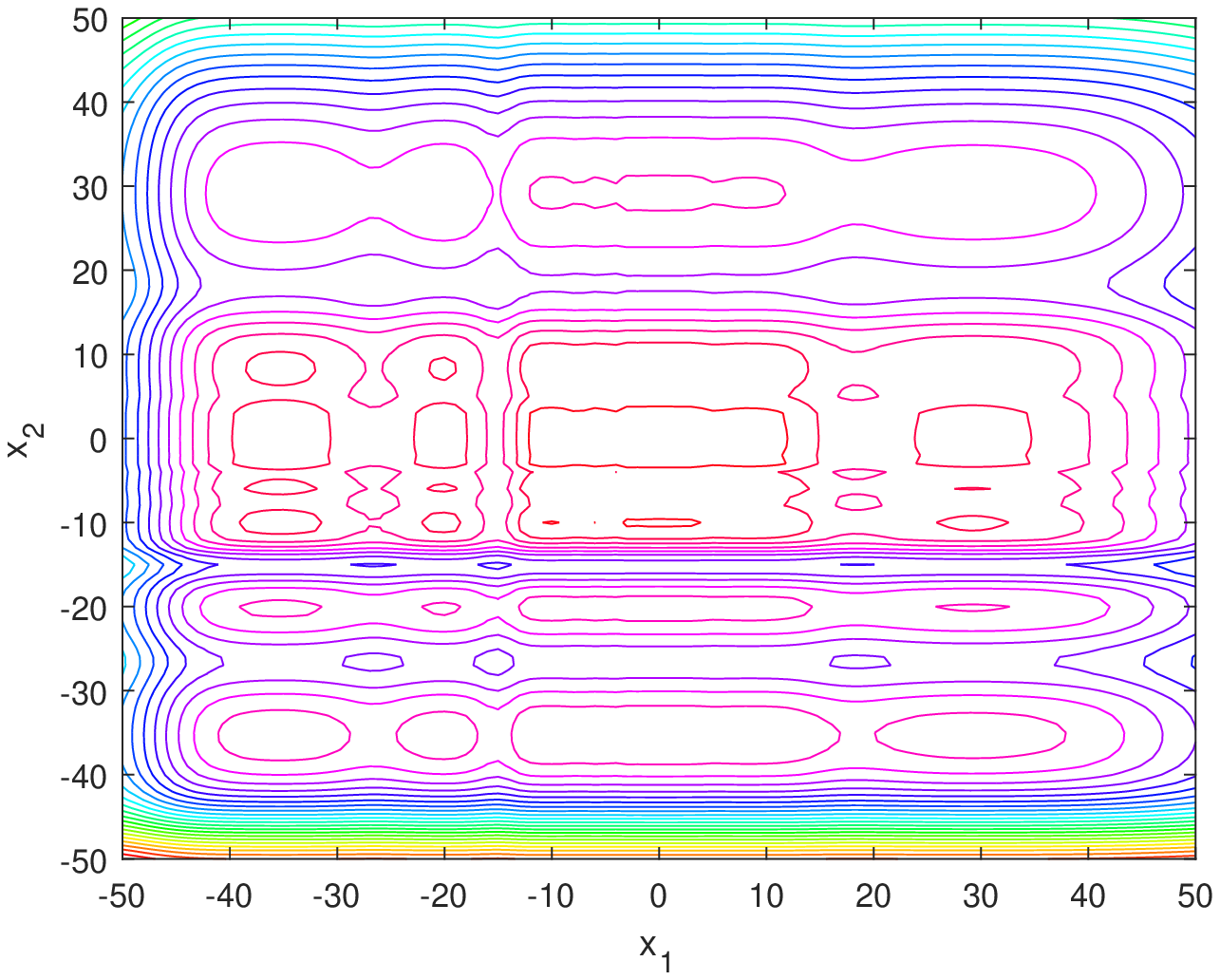}\label{fig:Component_irregular-assymetric-w73-contour}}
        \\
         \subfigure[{\scriptsize $\mathbf{w}=[10,2]$ (ill-conditioned).}]{\includegraphics[width=0.45\linewidth]{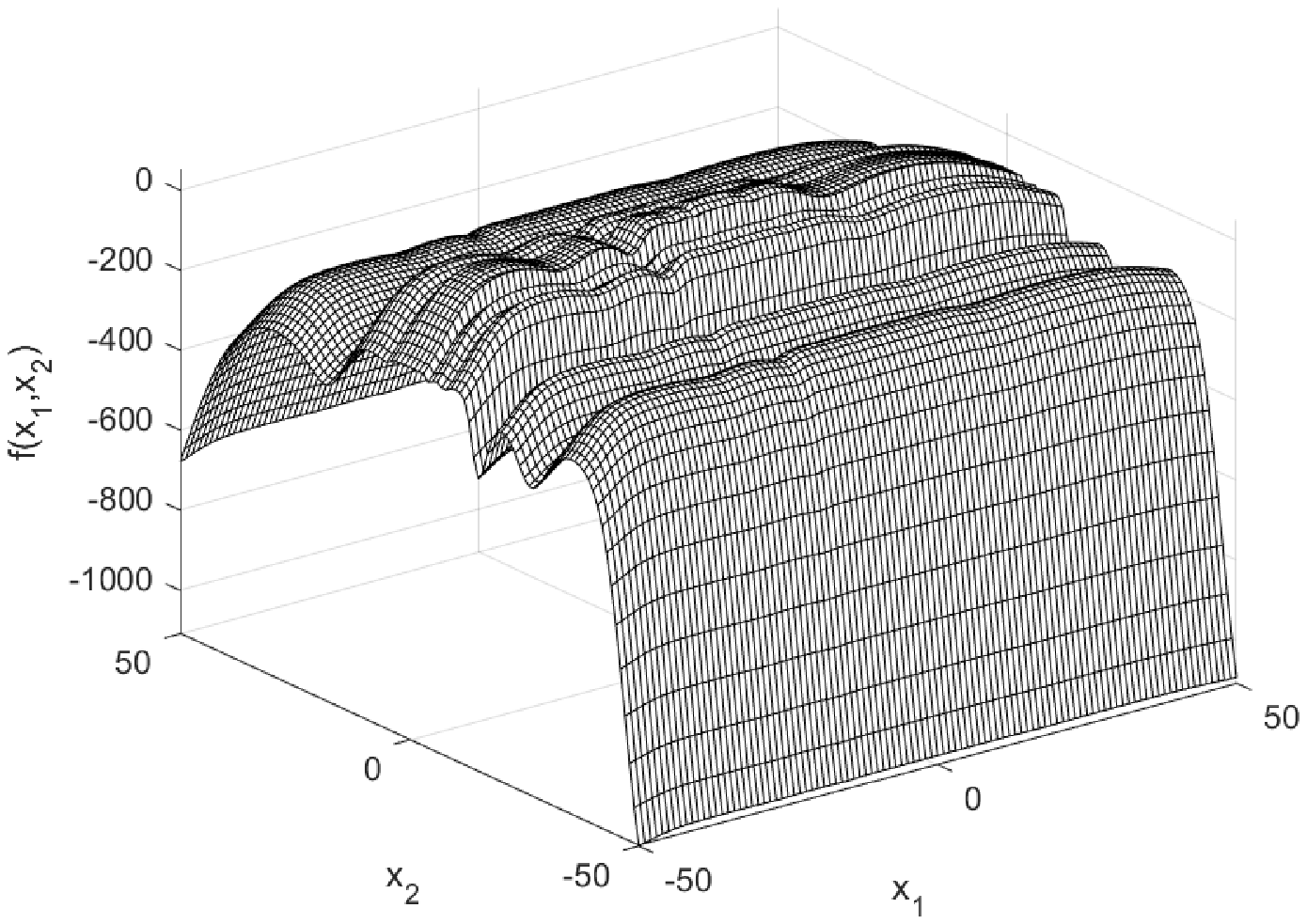}\label{fig:Component_irregular-assymetric-w102-surf}}
&
    \subfigure[{\scriptsize }]{\includegraphics[width=0.45\linewidth]{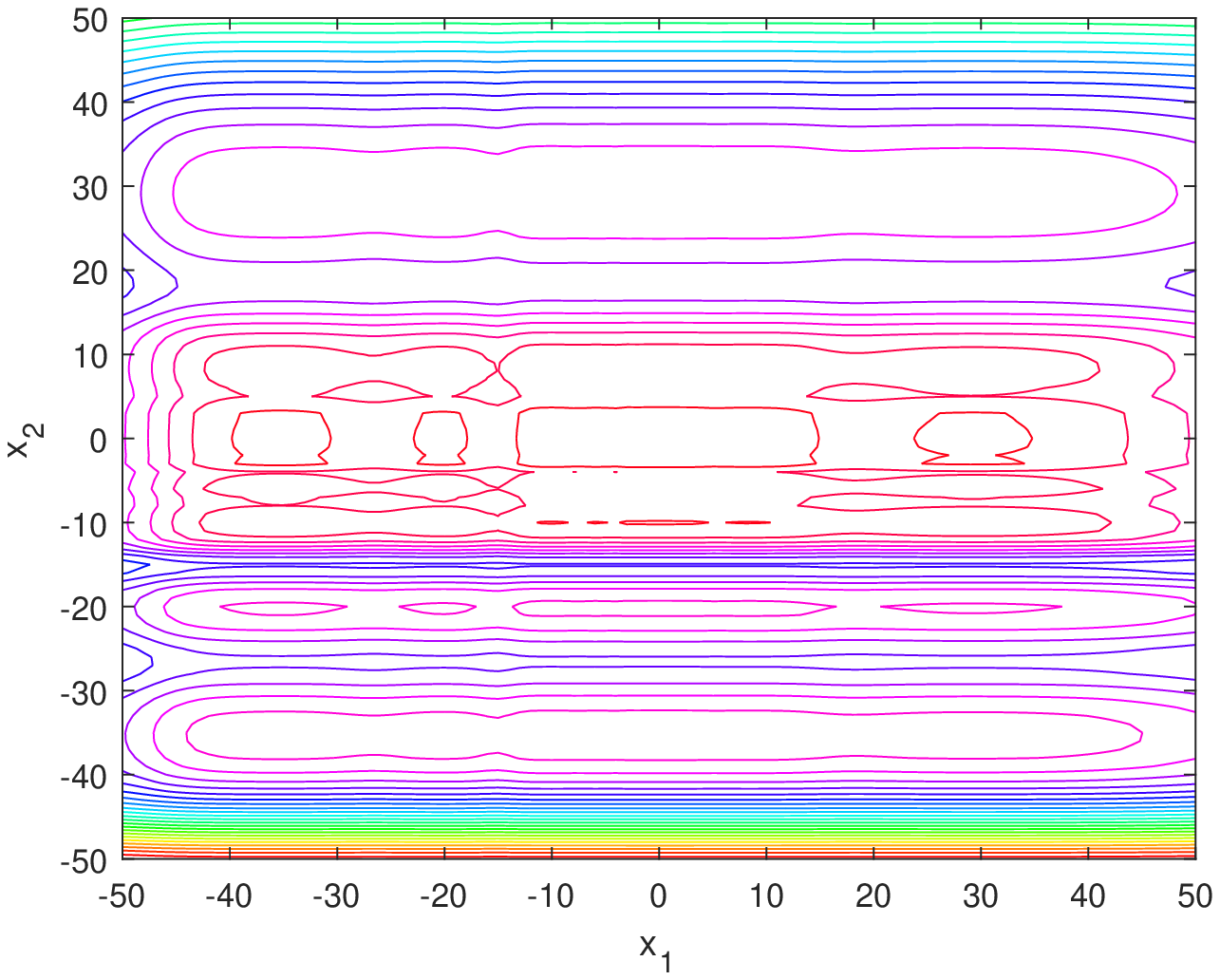}\label{Component_irregular-assymetric-w102-contour}}
\end{tabular}
\caption{Three components generated by \eqref{eq:irGMPB} whose $\mathbf{R}=\mathbf{I}$, $\tau=0.5$, and $\eta_{1,2,3,4}=[0,5,5,10]$.
The width values of these components are different. }
\label{fig:Cmponent:ill-conditioning}
\end{figure}

GMPB is capable of generating components with various condition numbers.
A component generated by GMPB has a width value in each dimension.
When the width values of a component are identical in all dimensions, the condition number of the component will be one, i.e., it is not ill-conditioned.
The condition number of a component is the ratio of its largest  width value to its smallest value~\cite{yazdani2020benchmarking}. 
If a component's width value is stretched in one axis's direction more than the other axes, then, the component is ill-conditioned.
Figure~\ref{fig:Cmponent:ill-conditioning} depicts three components with different condition numbers.

\begin{figure}[tp!]
\centering
\begin{tabular}{cc}
     \subfigure[{\scriptsize Not rotated, i.e., $\mathbf{R}=\mathbf{I}$.}]{\includegraphics[width=0.45\linewidth]{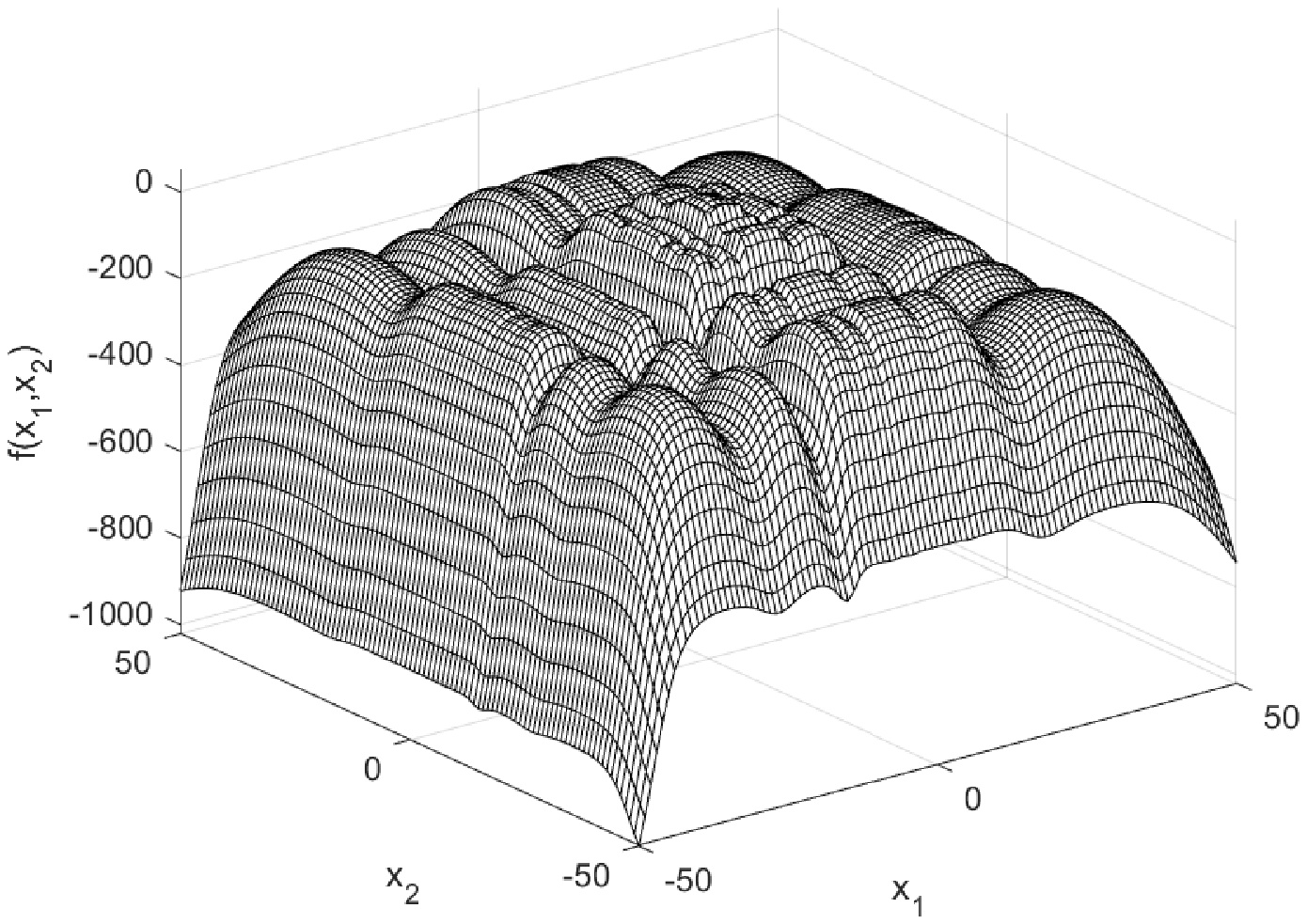}\label{fig:Component_irregular-assymetric-nonRotated-surf}}
&
    \subfigure[{\scriptsize }]{\includegraphics[width=0.45\linewidth]{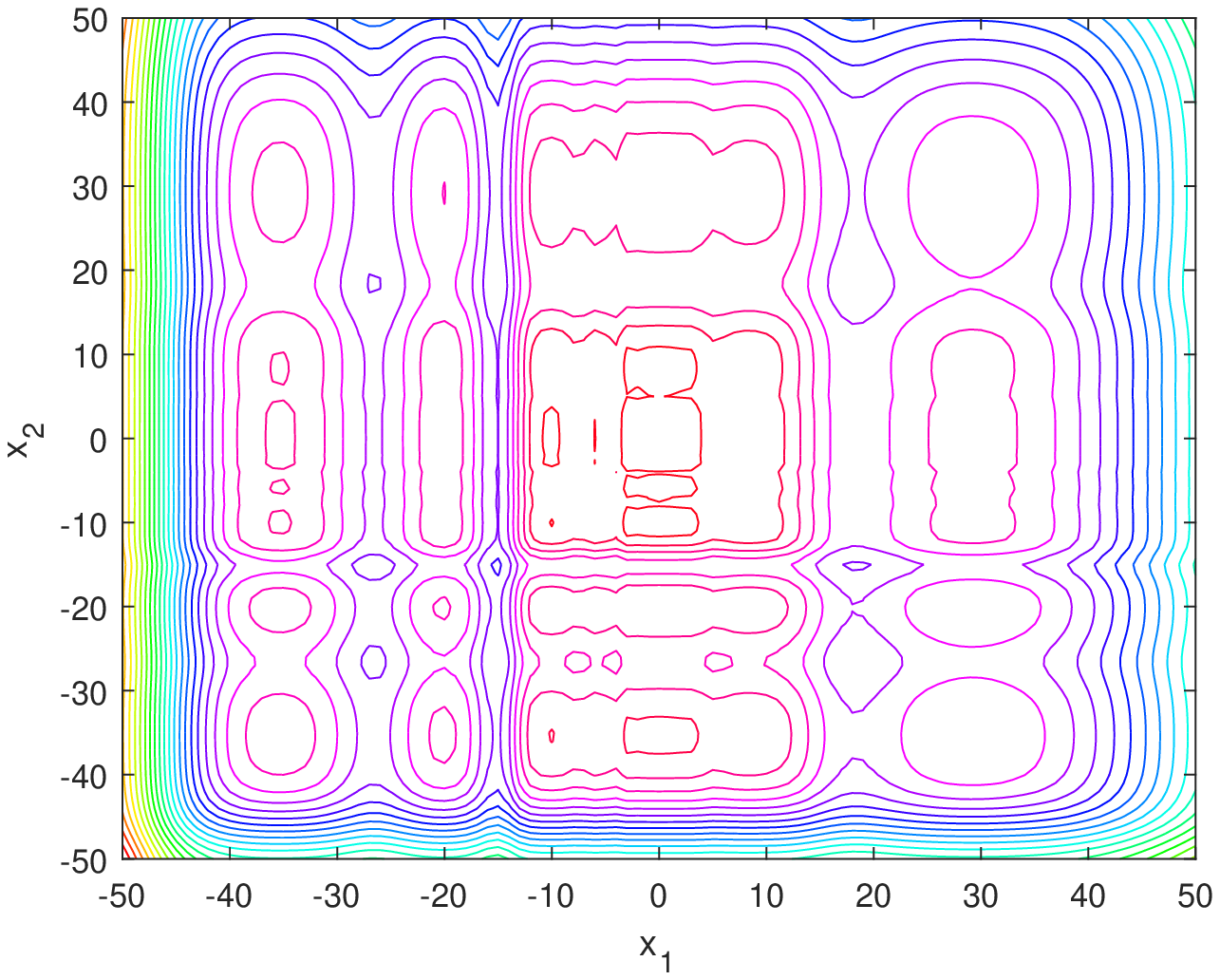}\label{fig:Component_irregular-assymetric-nonRotated-contour}}
    \\
         \subfigure[{\scriptsize Rotated ($45^\circ$)}]{\includegraphics[width=0.45\linewidth]{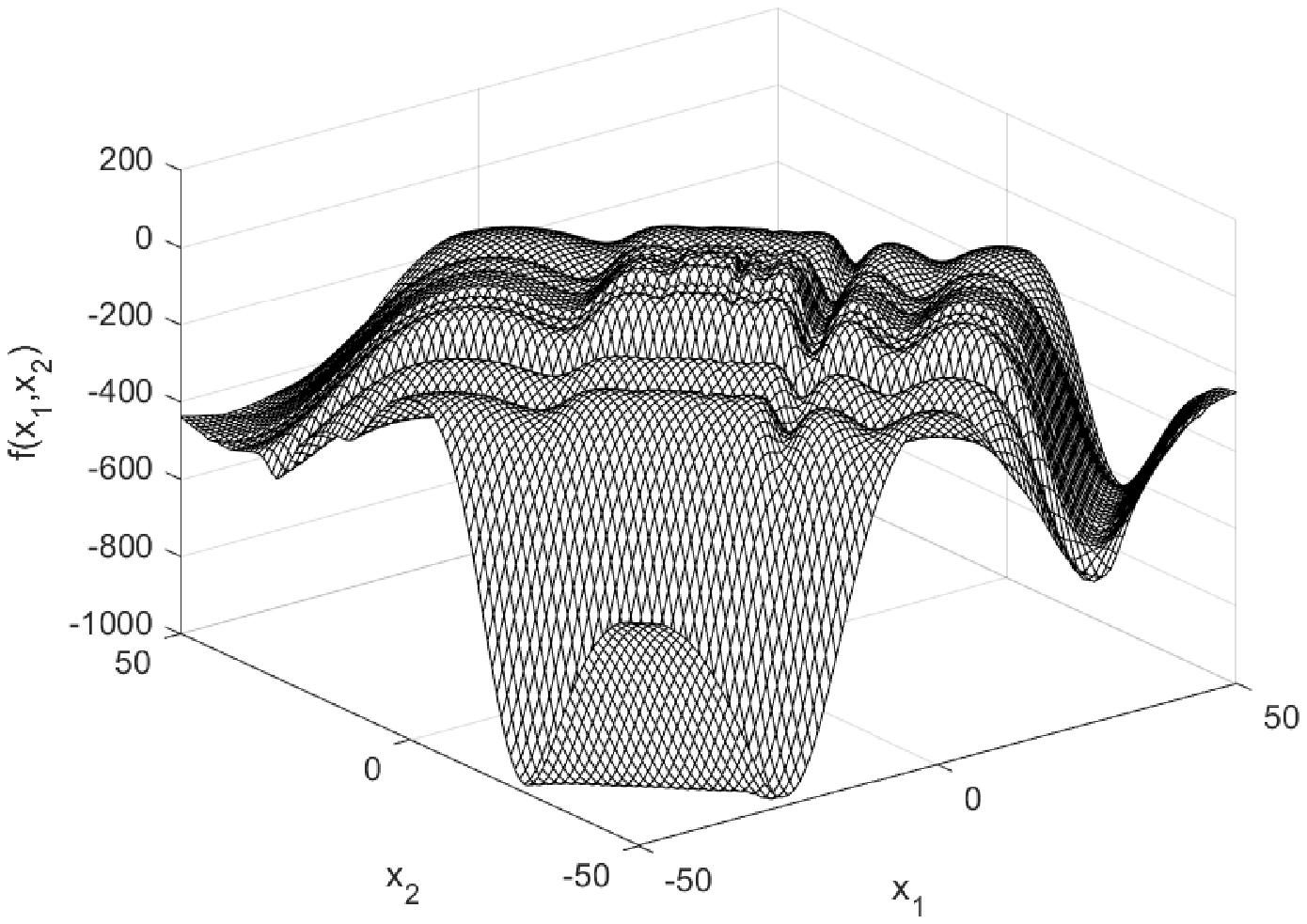}\label{fig:Component_irregular-assymetric-Rotated-surf}}
&
    \subfigure[{\scriptsize }]{\includegraphics[width=0.45\linewidth]{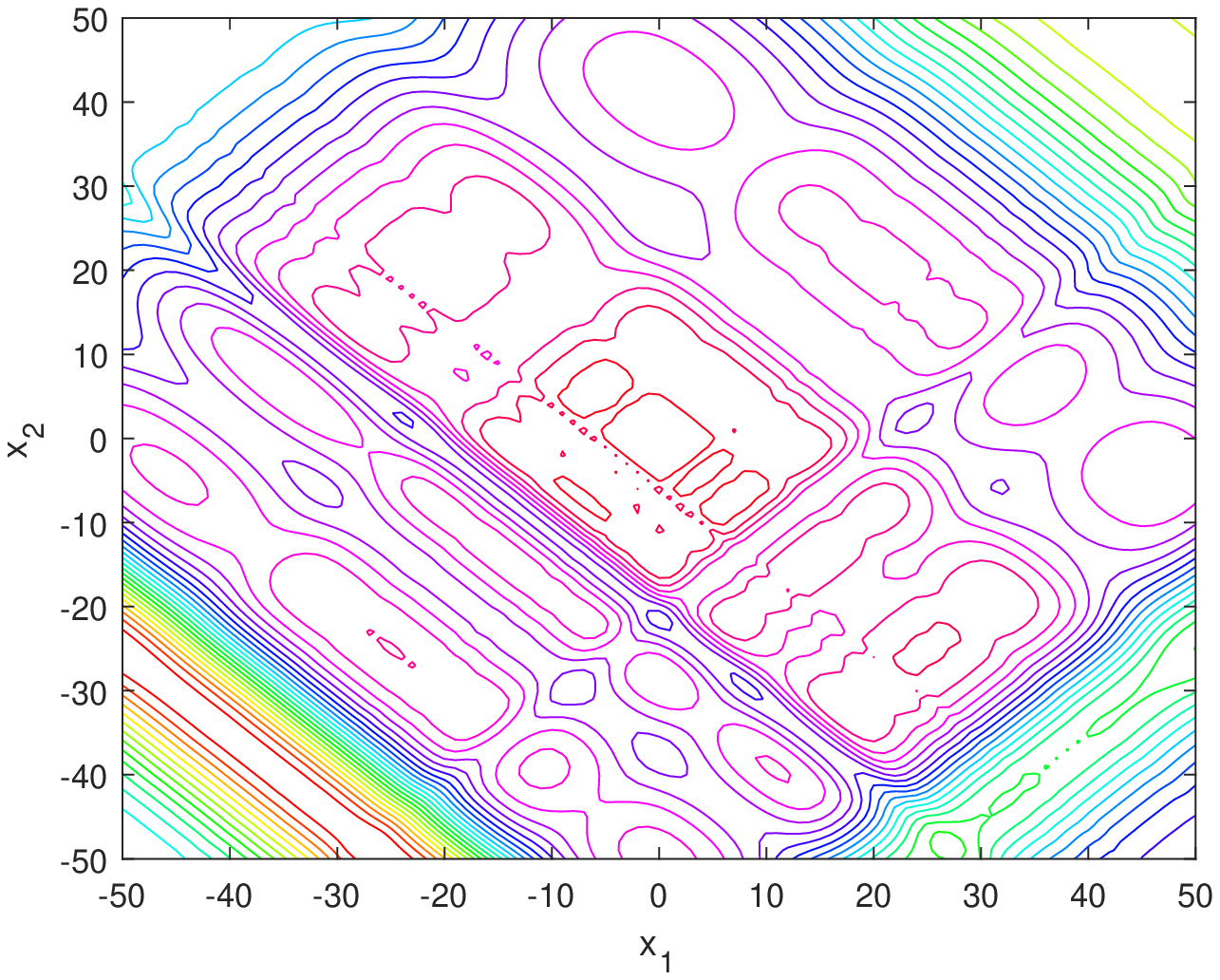}\label{fig:Component_irregular-assymetric-Rotated-contour}}
\end{tabular}
\caption{Two components generated by \eqref{eq:irGMPB} with different  $\mathbf{R}$ configurations.
The rest of parameters' values are identical for both components. }
\label{fig:Cmponent:rotation}
\end{figure}

In GMPB, each component $i$ is rotated using $\mathbf{R}_i$.
If $\mathbf{R}_i=\mathbf{I}$, then the component $i$ is not rotated (See Figure~\ref{fig:Component_irregular-assymetric-nonRotated-surf}).
Figure~\ref{fig:Component_irregular-assymetric-Rotated-surf} shows the component from Figure~\ref{fig:Component_irregular-assymetric-nonRotated-surf} which has been rotated by $\theta_i=\frac{\pi}{4}$.
In GMPB, by changing $\theta_i$ over time using \eqref{eq:angle}, the variable interaction degrees of the $i$th component change over time.

\subsubsection{Sub-function characteristics}

\begin{figure}[tp!]
\centering
\begin{tabular}{cc}
     \subfigure[{\scriptsize }]{\includegraphics[width=0.45\linewidth]{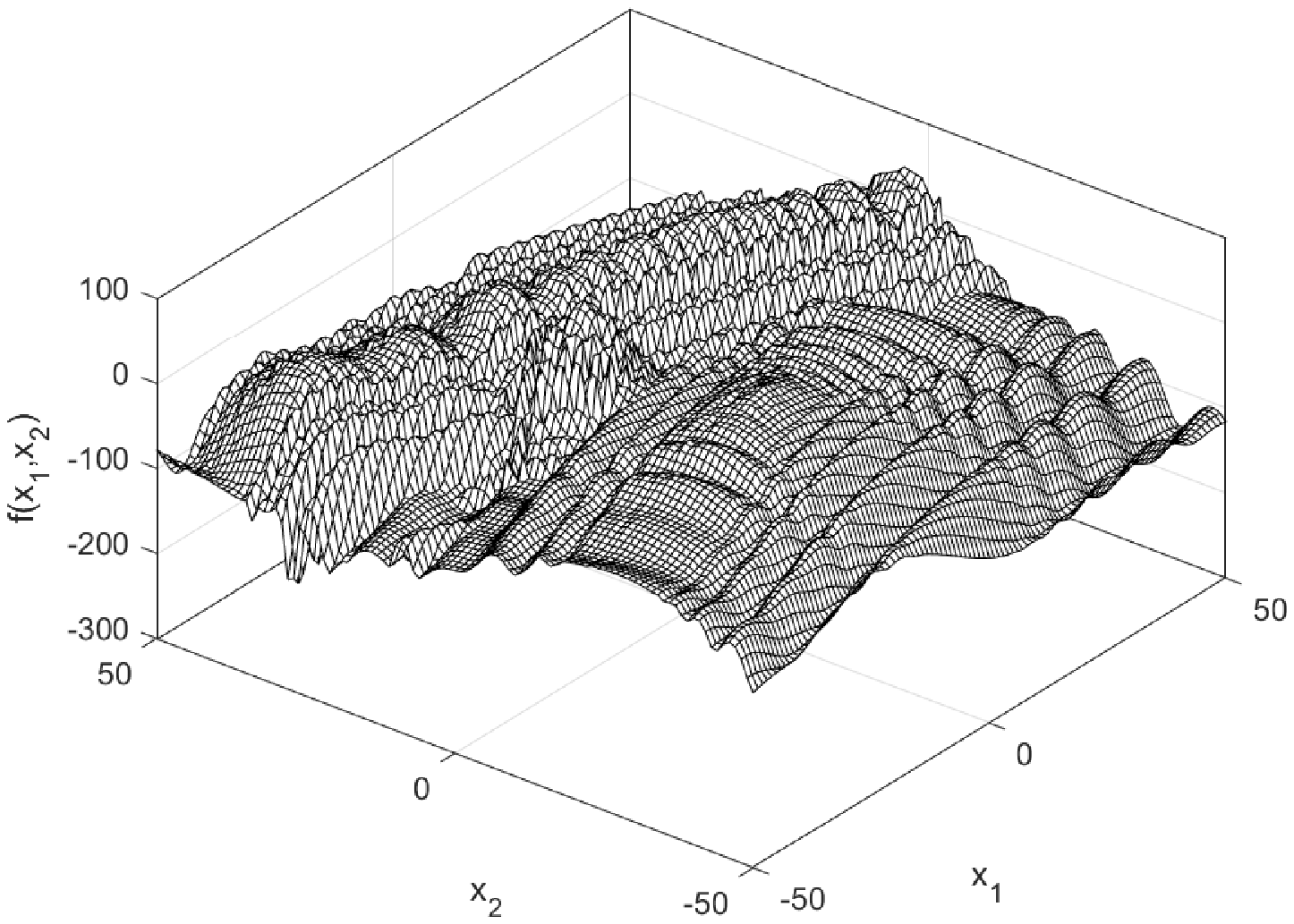}\label{fig:sub-function-surf}}
&
    \subfigure[{\scriptsize }]{\includegraphics[width=0.45\linewidth]{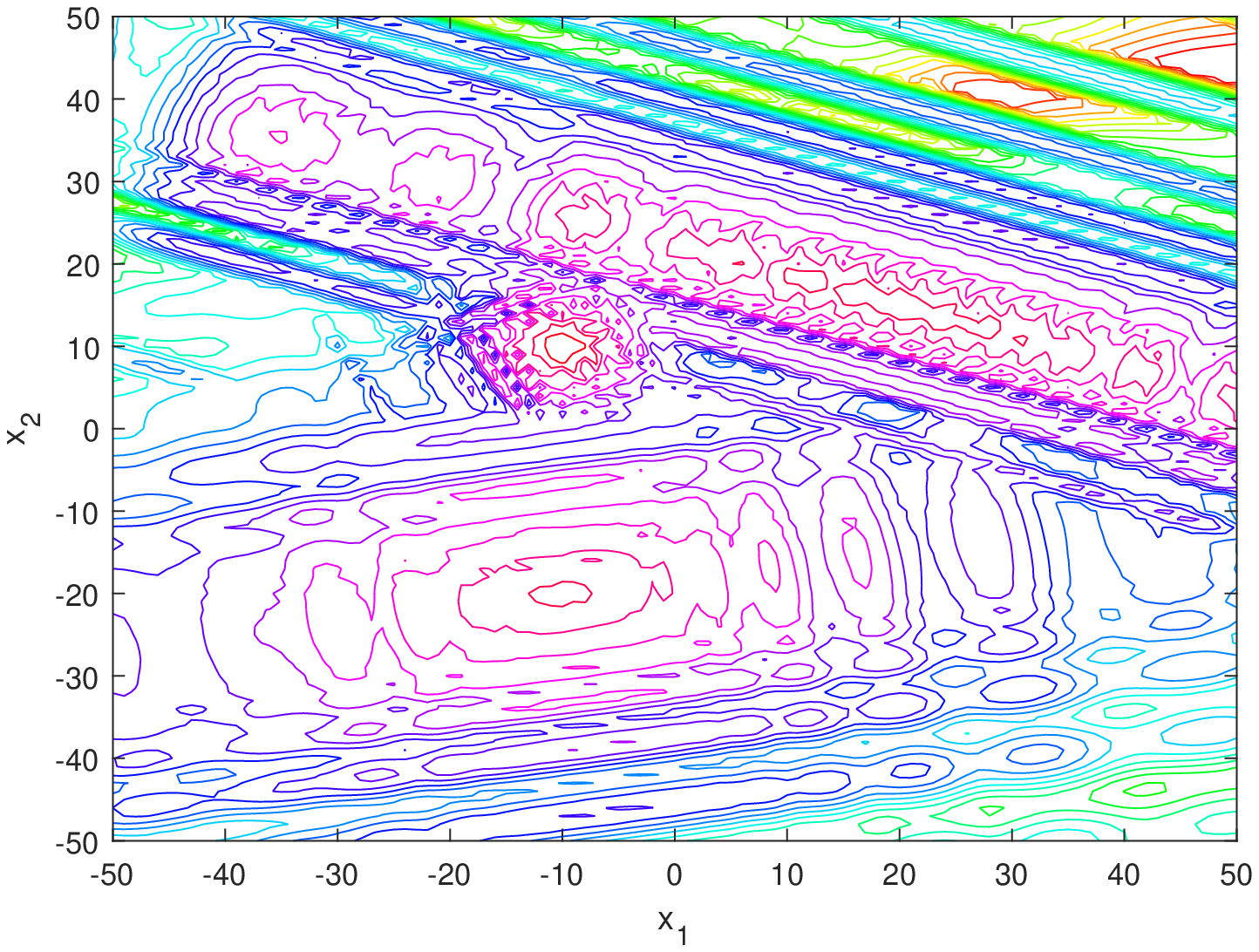}\label{fig:sub-function-contour}}
\end{tabular}
\caption{A sub-function generated by \eqref{eq:irGMPB} with three components whose parameter settings are different.}
\label{fig:sub-function}
\end{figure}

In GMPB, a sub-function is constructed by assembling several components using a $\max(\cdot)$ function in \eqref{eq:irGMPB}, which determines the basin of attraction of each component.
As shown in~\cite{yazdani2019scaling}, the landscape made by several components (i.e., $m>1$) by the $\max(\cdot)$ function is fully non-separable.
Figure~\ref{fig:sub-function} shows a landscape with three components.

\subsubsection{Problem characteristics}

A problem instance generated by GMPB is modular and constructed by assembling several sub-functions using \eqref{eq:wCMPB}.
These problems can be heterogeneous since the sub-functions can have different characteristics, such as different number of components, dimension, global optimum position, and change intensity.
Besides, by assigning different weight values to each sub-function $i$ ($\omega_i$), the problems can exhibit imbalance among its components.
In an imbalanced problem, the contributions of sub-functions on the overall fitness value are different. 
Consequently, some sub-functions become more important for the optimization algorithms~\cite{yazdani2018thesis}. 

\begin{figure}[tp!]
\centering
\begin{tabular}{cc}
     \subfigure[{\scriptsize 1-dimensional sub-function with two simple components ($\tau=0$, and $\eta_{1,2,3,4}=[0,0,0,0]$).}]{\includegraphics[width=0.45\linewidth]{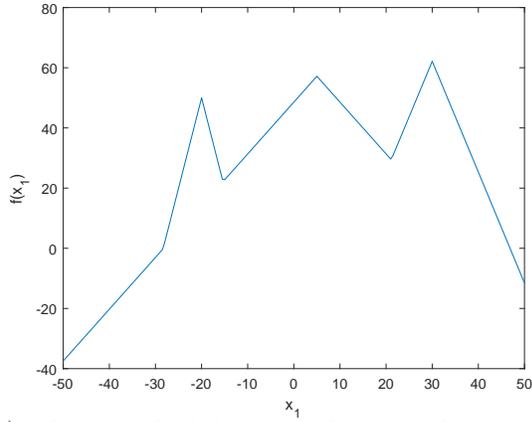}\label{fig:Problem1D_m3}}
&
    \subfigure[{\scriptsize 1-dimensional sub-function with three simple components ($\tau=0$, and $\eta_{1,2,3,4}=[0,0,0,0]$).}]{\includegraphics[width=0.45\linewidth]{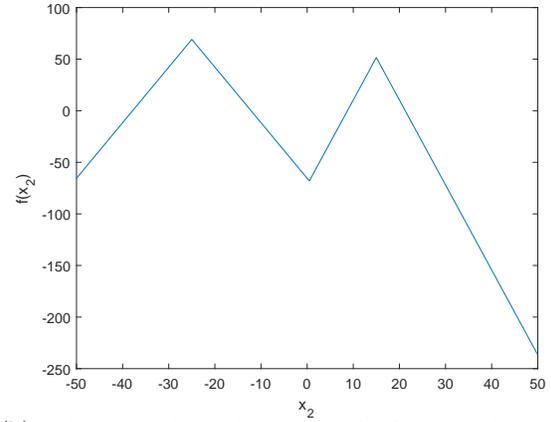}\label{fig:Problem1D_m2}}
    \\
         \subfigure[{\scriptsize The 2-dimensional landscape by assembling the landscapes shown in \ref{fig:Problem1D_m3} and \ref{fig:Problem1D_m2} by \eqref{eq:wCMPB}.}]{\includegraphics[width=0.45\linewidth]{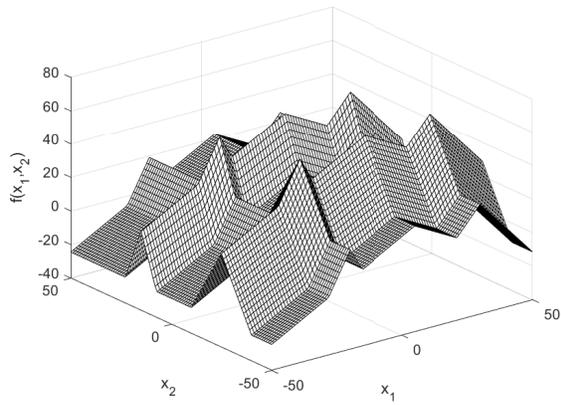}\label{fig:Problem-composed-surf}}
&
    \subfigure[{\scriptsize }]{\includegraphics[width=0.45\linewidth]{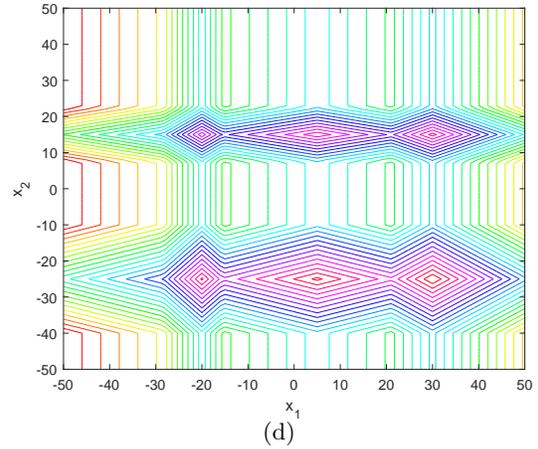}\label{fig:Problem-composed-contour}}

\end{tabular}
\caption{Exponentially growing the number of promising regions by assembling sub-functions using \eqref{eq:wCMPB}. }
\label{fig:GrowOfPeakNumbers}
\end{figure}

A consequence of GMPB's design is the exponential growth in the total number of promising regions that can contain the global optimum in a future environment. 
In each GMPB's sub-function, the total number of such promising regions can be up to the number of components whose center position can become the global optimum after environmental changes.
However, by assembling several sub-functions using \eqref{eq:wCMPB}, the number of such promising regions in the problem becomes:
\begin{align}
\label{eq:CMPBpeaknumber}
M = \prod_{i=1}^{n+l}m_i,
\end{align}
where $m_i$ is the number of components in the $i$th sub-function. 
It should be noted that $M$ is the maximum number of promising regions that can exist in the landscape, which may change over time due to coverage of smaller components by larger ones. 
For the sake of clarity, we provide an illustrative example. 
In Figure~\ref{fig:Problem1D_m3} and Figure~\ref{fig:Problem1D_m2}, two 1-dimensional sub-functions with two and three components (in the simplest form, i.e., conical peaks) are shown. 
The 2-dimensional function constructed based on~\eqref{eq:wCMPB} with $\omega_1=\omega_2=1$ results in a total of $2\times3=6$ promising regions.

\section{Problem instances}
\label{sec:senario}

We present 15 GMPB large-scale scenarios with five different variable interaction structures in 50-, 100-, and 200-dimensional spaces, which are shown in Table~\ref{tab:scenarios}. 
Note that in dynamic environments, the curse of dimensionality happens in much lower dimensions in comparison to  static optimization problems~\cite{yazdani2019scaling}.
In fact, the very limited available computational resources (i.e., the number of fitness evaluations) in each environment of DOPs results in a significant deterioration in the performance of locating and tracking the moving global optimum even in problems with 50-dimensional search space.

In the presented 15 GMPB large-scale scenarios in Table~\ref{tab:scenarios}, the number of dimensions of the main function $F^{(t)}$ ($d$) in \eqref{eq:wCMPB} and all sub-functions $f^{(t)}_i$  ($d_i$), and also the variable interaction structures of the main function $F^{(t)}$  are fixed and do not change over time. The order of components in the actual benchmarks is taken from a random permutation; however, in Table~\ref{tab:scenarios} we present them in sorted order for better readability.
Every $\vartheta$ fitness evaluations, the center position, height, width vector, angle, and irregularity parameters of each component in each sub-function change using the dynamics presented in~\eqref{eq:center}~to~\eqref{eq:eta}.
The parameter settings of each sub-function $f_i^{(t)}$  \eqref{eq:irGMPB} including  their change severity values (in~\eqref{eq:center}~to~\eqref{eq:eta}) and their ranges are listed in Table~\ref{tab:parameterSettings-part1}.
Besides, the temporal parameter settings are shown in Table~\ref{tab:parameterSettings-part2}.
Finally, Table~\ref{tab:InitialparameterValues} shows the initial values of the components' parameters.

The problem instances can become more difficult by setting some parameters to what we call \emph{challenging settings}. 
These challenging settings are characterized as follows.
By increasing the shift severity values, optima relocate more severely and tracking them will become more time-consuming.
Increasing the number of components results in a larger number of promising regions that may contain the global optimum after environmental changes.
Therefore, dynamic optimization algorithms (DOAs) that try to locate and track multiple moving promising regions~\cite{yazdani2020adaptive} will face difficulties in tackling such problems. 
In fact, covering larger numbers of promising regions is more computational resource consuming which is challenging due to the limited available computational resources in each environment.
Another factor to increase the difficulty of each sub-function is to set all irregularity parameters' boundaries to the challenging settings.
This parameter setting results in larger plateau (See Figure~\ref{fig:Cmponent:irregular}), which significantly increases the probability of premature convergence.
Finally, by decreasing the value of $\vartheta$, we can increase the change frequency (see Table~\ref{tab:parameterSettings-part2}).
In problem instances with higher change frequencies, the available computational resources in each environment (i.e., the number of fitness evaluations in each environment) are more limited which results in increased difficulty.

\begin{table}[t]
\small
\centering
\caption{The variable interaction structures and dimensionality of the 15 GMPB Scenarios. 
The order of components in the actual benchmarks is taken from a random permutation. 
Here we present them in sorted order for better readability.}
\label{tab:scenarios}
\begin{tabu}{@{ }c@{ }c@{\ \ }l@{ }c@{ }}
\toprule
Function & $d$ & Dimensionality of Nonseparable Components & \# Separable Vars \\ \midrule
$f_1$    & 50  & $\{2,3,5,6,7,8,10\}$                      & 10        \\
$f_2$    & 50  & $\{2,3,5,5\}$                             & 35        \\
$f_3$    & 50  & $\{2,2,3,5,5,5,5,5,8,10\}$                & 0         \\
$f_4$    & 50  & ---                                       & 50        \\
$f_{5}$ & 50  & $\{50\}$                                   & 0         \\ 
\midrule
$f_{6}$ & 100 & $\{2,2,3,5,5,6,6,8,8,10,10,15\}$           & 20        \\
$f_{7}$ & 100 & $\{2,2,3,3,5,5,10\}$                       & 70        \\
$f_{8}$ & 100 & $\{2,2,2,2,3,3,5,5,5,5,5,5,8,8,10,10,20\}$ & 0         \\
$f_{9}$ & 100 & ---                                        & 100       \\
$f_{10}$ & 100 & $\{100\}$                                 & 0         \\
\midrule
$f_{11}$ & 200 & $\{2 , 2 , 3 , 5 , 5 ,  6 , 6 , 8 , 8 , 10 , 10 , 15 , 20 , 20 , 30\}$           & 50        \\
$f_{12}$ & 200 & $\{2 , 3 , 5 , 10 , 20 , 30\}$                       & 130        \\
$f_{13}$ & 200 & $\{2 , 2 , 2 , 3 , 5 , 5 , 5 , 5 , 5 , 8 , 8 , 10 , 10 , 10 , 20 , 20 , 30 , 50\}$ & 0         \\
$f_{14}$ & 200 & ---                                        & 200       \\
$f_{15}$ & 200 & $\{200\}$                                  & 0         \\
\bottomrule
\end{tabu}
\end{table}

\begin{table*}[tp!] 
    \small
\centering
  \caption{Parameter settings of sub-functions generated by \eqref{eq:irGMPB} and their dynamics. }
  \label{tab:parameterSettings-part1}
  \begin{tabular}{ll@{\ \ }*{2}{c}}
    \toprule
    Parameter & Symbol & Default setting & Challenging setting\\ 
\midrule
Shift severity   & $\tilde{s}_i$	&  $ \mathcal{U}[1,3]$   &   $ \mathcal{U}[3,5]$\\
Numbers of components  & $m_i$	             & $ \mathcal{U}[5,15]$  &   $ \mathcal{U}[15,35]$\\
Angle severity     & $\tilde{\theta}_i$ & $ \mathcal{U}[\pi/12,\pi/6]$ & -\\
Height severity   & $\tilde{h}_i$	    & $\mathcal{U}[5,9]$& - \\
Width severity    & $\tilde{w}_i$	    & $\mathcal{U}[0.5,1.5]$& -\\
Irregularity parameter $\tau$ severity    & $\tilde{\tau}_i$	    &  $ \mathcal{U}[0.05,0.15]$ & - \\
Irregularity parameter $\eta$ severity    & $\tilde{\eta}_i$	    &  $ \mathcal{U}[1,3]$ & - \\
Weight of sub-function $i$       &	$\omega_i$   & $ \mathcal{U}[0.5,3]$ & - \\
Search range                        &   $[Lb,Ub]^{d_i}$	    &    $[-50,50]^{d_i}$ & - \\
Height range                         &  $[h_{\mathrm{min}},h_\mathrm{max}]$  	    &    $[30,70]$ & -  \\
Width range                          &   $[w_\mathrm{min},w_\mathrm{max}]^{d_i}$  	    &   $[1,12]^{d_i}$ & - \\
Angle range                           & $[\theta_\mathrm{min},\theta_\mathrm{max}]$ & $[-\pi,\pi]$  & - \\
Irregularity parameter $\tau$ range   &   $[\tau_\mathrm{min},\tau_\mathrm{max}]$   	    & $[-0.5,0.5]$  & - \\
Irregularity parameter $\eta$ range   &  $[\eta_\mathrm{min},\eta_\mathrm{max}]$   	    & $[-20,20]$  & -  \\
    \bottomrule
  \end{tabular}
 \end{table*}

\setlength\tabcolsep{4pt} 
\begin{table}[tp!] 
    \small
\centering
  \caption{Temporal parameter settings. }
  \label{tab:parameterSettings-part2}
  \begin{tabular}{llcc}
    \toprule
    Parameter &Symbol & Default setting & Challenging setting\\ 
    \midrule
Change frequency   & $\vartheta$	             & $500 \cdot d$ & $200 \cdot d$\\
Number of Environments                 &   $T$       &  30  & - \\
    \bottomrule
  \end{tabular}
\end{table}
\setlength\tabcolsep{6pt} 

\setlength\tabcolsep{4pt} 
\begin{table}[tp!] 
    \small
\centering
  \caption{Initial values of components' parameters. }
  \label{tab:InitialparameterValues}
 \begin{threeparttable}
  \begin{tabular}{llcc}
    \toprule
    Parameter & Symbol & Initial value\\ 
    \midrule
Center position   & $\vec{c}^{(0)}_i$	             &  $\mathcal{U}[Lb,Ub]^{d_i}$\\
Height                 &   $h^{(0)}_i$       &  $\mathcal{U}[h_{\mathrm{min}},h_\mathrm{max}]$ \\
Width                          & $\vec{w}^{(0)}_i$ & $\mathcal{U}[w_\mathrm{min},w_\mathrm{max}]^{d_i}$  \\
Angle                            & $\theta^{(0)}_i$ & $\mathcal{U}[\theta_\mathrm{min},\theta_\mathrm{max}]$ \\
Irregularity parameter $\tau$    & $\tau^{(0)}_i$    &  $\mathcal{U}[\tau_\mathrm{min},\tau_\mathrm{max}]$  \\
Irregularity parameter $\eta$    &$\eta^{(0)}_i$    &  $\mathcal{U}[\eta_\mathrm{min},\eta_\mathrm{max}]$   \\    
Rotation matrix    &$\mathbf{R}^{(0)}_i$    &   $\mathrm{GS}(\mathcal{N}(0,1)^{d_i\times d_i})$\tnote{\textbf{$\dagger$}}  \\   
 \bottomrule
  \end{tabular}
   \begin{tablenotes}
   \begin{scriptsize}
  \item[$\dagger$] $\mathbf{R}_i$ is initialized by performing the Gram-Schmidt orthogonalization method $\mathrm{GS}(\cdot)$ on a $d_i\times d_i$ matrix with normally distributed entries.
   \end{scriptsize}
    \end{tablenotes}
 \end{threeparttable}
\end{table}
\setlength\tabcolsep{6pt} 

\section{Performance indicator}
\label{sec:PerformanceIndicator}

To measure the performance of DOAs in solving the problem instances generated by GMPB, the average error of the best found solutions at the end of all environments (i.e., the best before change error ($E_\mathrm{BBC}$)) is used as the performance indicator~\cite{trojanowski1999searching}: 
\begin{align}
\label{eq:performance}
E_\mathrm{BBC}= \frac{1}{T} \sum_{t=1}^T\left(d^{-1}\sum_{i=1}^n \omega_i d_i h_{i,\max}^{(t)}-F^{(t)}\left(\mathbf{g}^{*(t)}\right)\right),
\end{align}
where $\mathbf{g}^{*(t)}$ is the best found position in the $t$th environment which is fetched at the end of the environment, $F^{(t)}$ is \eqref{eq:wCMPB}, $h_{i,\max}^{(t)}$ is the maximum height value among the components of the $i$th sub-function in the $t$th environment, and $d^{-1}\sum_{i=1}^n \omega_i d_i h_{i,\max}^{(t)}$ calculates the global optimum fitness value in the $t$th environment.

\section{Source code}

The MATLAB\footnote{Version R2019a} source code of the large-scale scenarios generated by the GMPB can be downloaded from~\cite{yazdani2021PSOCTR4LSGMPB}.
As a sample optimization algorithm, we have added $\mathrm{PSO}_\mathrm{CTR}$ to this code as the optimizer.
$\mathrm{PSO}_\mathrm{CTR}$ is a cooperative coevolutionary algorithm, which has been designed to tackle large-scale DOPs~\cite{yazdani2019scaling}.
This algorithm uses the DG2~\cite{omidvar2017dg2} for decomposing the problem and identifying the variable interaction structure of the problem.
Thereafter, each sub-function is assigned to a multi-population PSO algorithm which work cooperatively via a \emph{context vector}~\cite{vandenBergh2004cooperative}.

\bibliography{bib} 
\bibliographystyle{IEEEtran}

\end{document}